\documentclass{article} 
\usepackage{texdraw}\usepackage{graphicx}
\usepackage{amsmath} 
\usepackage{color}
\usepackage{colordvi}
\usepackage{texdraw}\usepackage{graphicx}
\usepackage{dcpic,pictex}

\definecolor{purple}{cmyk}{0.45,0.86,0,0}\definecolor{brickred}{cmyk}{0,0.89,0.94,0.28}
\definecolor{maroon}{cmyk}{0,0.87,0.68,0.32}\definecolor{lyellow}{cmyk}{0,0,0.68,0}
\definecolor{magenta}{cmyk}{0,1,0,0}
\definecolor{zals}{cmyk}{0.75,0.12,0.98,0.5}
\definecolor{pamats}{cmyk}{0.54,0.00,0.84,0.0}
\definecolor{brickred}{cmyk}{0,0.89,0.94,0.28}
\definecolor{cadetblue}{cmyk}{0.62,0.57,0.23,0}
\definecolor{teo}{cmyk}{0.63,0.83,0.33,0.13}

\def\newpic#1{%
   \def\emline##1##2##3##4##5##6{%
      \put(##1,##2){\special{em:point #1##3}}%
      \put(##4,##5){\special{em:point #1##6}}%
      \special{em:line #1##3,#1##6}}}
\newpic{}

\usepackage[cp1251,cp1257]{inputenc}\usepackage[english,latvian]{babel}\usepackage{amssymb,latexsym,amsfonts}
\newcommand{\apz}{\protect\makebox[1.8\width]{\rule[3.3pt]{8.5pt}{0.4pt} \hspace*{-12.6pt}\protect\raisebox{-.3\height}{$\leftharpoondown$}}}

\newtheorem{teor}{Theorem}[section]

\newtheorem{corollary}[teor]{Corollary}
\newtheorem{lemma}[teor]{Lemma}
\newtheorem{proposition}[teor]{Proposition}
\newtheorem{definition}[teor]{Definition}

\newtheorem{example}[teor]{Example}

\input cyracc.def

\begin{document}
\selectlanguage{english}
\title{Machine $B_4$}
\date{}
\author{J\=anis Buls \\
{\small Department  of  Mathematics, University of Latvia, Jelgavas iela 3,}\\
{\small R\=\i ga, LV-1004 Latvia,
buls@edu.lu.lv
}}
\def\keywords{\begin{center}{\bf Keywords}\end{center}
{automata (machines) groups, dense orbit, topological transitivity}}

\maketitle
\title{}
\abstract{\small{We construct map $\xi$. It exhibits dense orbits for all $x\in\overline{0,1}^\omega$. We give elementary proofs for all statements. }}


\section{ Preliminaries} 
\normalsize Let $A$ be a finite non-empty set and $A^*$ 
the free monoid generated by $A$. The set $A$ is also called an {\em 
alphabet}, its  elements are called {\em letters} and 
those of $A^*$ are called {\em finite words}. The identity element of $A^*$ is called 
an {\em empty word} and denoted by $\lambda$. We set 
$A^+=A^*\backslash\{\lambda\}$.

A word $w\in A^+$ can be written uniquely as a sequence of letters as 
$w=w_1w_2\ldots w_l$, with $w_i\in A$, $1\le i\le l$, $l>0$. The integer 
$l$ is called the {\em length} of $w$ and denoted by $|w|$. The length of 
$\lambda$ is 0. We set $w^0=\lambda$ and $\forall i \in \mathbb{N} \; w^{i+1}=w^iw\,.$ 

The word $w'\in A^*$ is a {\em factor} (or {\em subword}) of 
$w\in A^*$ if there exists $u,v\in A^*$ such that $w=uw'v$.
The words $u$ and $v$ are called, respectively, a {\em prefix} and a {\em suffix}. A pair $(u,v)$ is called an {\em occurrence} of $w'$ in $w$.
A factor $w'$ is called {\em proper} if $w\ne w'$. We denote, respectively, by F$(w)$, Pref$(w)$ and Suff$(w)$ the sets of $w$ factors, prefixes and suffixes.

An (indexed) infinite word $x$ on the alphabet $A$ is any total mapping 
$x\,:\,\mathbb{N}\rightarrow A$. We shall set for any $i\ge0$, $x_i=x(i)$
and write
\[
x=(x_i)=x_0x_1\ldots x_n\ldots \;.
\]
The set of all the infinite words over $A$ is denoted by $A^\omega$.

The word $w'\in A^*$ is a {\em factor} of $x\in A^\omega$ if there exists 
$u\in A^*$, $y\in A^\omega$ such that $x=uw'y$. 
The words $u$ and $y$ are called, respectively, a {\em prefix} and a {\em suffix}.
We denote, respectively, by F$(x)$, Pref$(x)$ and Suff$(x)$ the sets of $x$ factors, prefixes and suffixes.
We write $u\smallsetminus x$ if $u\in{\rm F}(x)$.  For any $0\le m\le n$, $x[m,n]$ denotes a factor $x_mx_{m+1}\ldots x_n$. The word 
$x[m,n]$ is called an {\em occurrence} of $w'$ in $x$ if $w'=x[m,n]$. The 
suffix $x_nx_{n+1}\ldots x_{n+i}\ldots$ is denoted by $x[n,\infty)$.

If $v\in A^+$, then we denote by $v^\omega$ the infinite word
\[
v^\omega=vv\ldots v\ldots \;.
\]

The {\em concatenation} of $u=u_1u_2\ldots u_k\in A^*$ and $x\in 
A^\omega$ is the infinite word 
\[
ux=u_1u_2\ldots u_kx_0x_1\ldots x_n\ldots
\]
For denoting concatenation we sometimes use symbol \#. 

We use notation $\overline{0,n}$ to denote the set $\left \{ 0,1,...,n \right \}$. 
\section{Machine $B_4$}\label{n5.2}

\begin{figure}[t]
\unitlength 1mm
\begin{picture}(100,30)(10,82)
\put(60,105){\line(0,-1){20}}
\put(60,85){\line(1,0){25}}
\put(85,85){\line(0,1){20}}
\put(85,105){\line(-1,0){25}}
\put(60,95){\line(-1,0){25}}
\put(100,95){\line(-1,0){15}}
\put(45.00,95.00){\vector(-1,0){0.2}}
\put(90.00,95.00){\vector(-1,0){0.2}}
\put(72.67,95.00){\makebox(0,0)[cc]{$\langle Q,A,B \rangle$}}
\put(100.00,100.00){\makebox(0,0)[cc]{$u_1u_2\ldots u_n$}}
\put(47.00,100.00){\makebox(0,0)[cc]{$v_1v_2\ldots v_n$}}
\end{picture}
\caption{An abstract Mealy machine. }
\label{z51}
\end{figure}
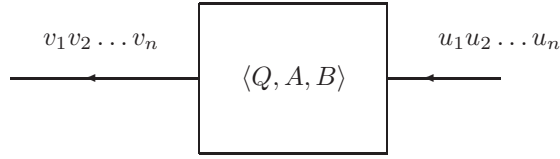

\begin{definition}
A 3-sorted algebra $V=\langle Q,A,B,\circ ,\ast \rangle$ is called  a 
{\em Mealy machine} if $Q,A,B$ are finite, nonempty sets, the mapping 
$Q\times A \stackrel{\circ}{\longrightarrow}Q$
is a total function and the mapping 
$Q\times A \stackrel{\ast }{\longrightarrow}B$
is a total surjective function. 
\end{definition}

If $A=B$ we do not insist on surjectivity of the map $\ast$. 
The, set $Q$ is called {\em state set}, sets $A,B$ are called {\em input} and {\em output alphabet}, respectively.
The mappings $\circ$ and $\ast$ may be extended to $Q\times A^*$  by defining
\[
\begin{array}{lr}
q\circ \lambda =q,\quad & q\circ (ua)=(q\circ u)\circ a, \\
q\ast \lambda  =\lambda ,\quad & q\ast (ua)=(q\ast u)\#((q\circ u)\ast a)\,,
\end{array}
\] 
for each $q\in Q$, $(u,a)\in A^*\times A$. See \ref{z51}\,fig. for interpretation of Mealy machine as a word transducer. 
Henceforth, we shall omit parentheses if there is no danger of confusion. So, for example, we will write $q\circ u\ast a$ instead of $(q\circ u)\ast a.$ Similarly, we will write $q\circ q^{'}\ast a$ instead of $q\circ (q^{'}\ast a)$ where $q^{'} \in Q$.

\begin{figure}[h]
\special{em:linewidth 0.4pt}
\unitlength 1.00mm
\linethickness{0.4pt}
\begin{picture}(131,58)(10,88)
\put(40.00,100.00){\circle{6.00}}
\put(100.00,130.00){\circle{6.00}}
\put(100.00,100.00){\circle{6.00}}
\put(40.00,130.00){\circle{6.00}}

\qbezier(42,132)(70,142)(98,132.00)
\put(50.00,134.45){\vector(-4,-1){0.2}}
\qbezier(42,128)(70,118)(98,128.00)
\put(90.00,125.55){\vector(4,1){0.2}}

\put(43,100){\line(1,0){54}}
\put(65.00,100){\vector(1,0){0.2}}

\put(40,127){\line(0,-1){24}}
\put(40.00,110){\vector(0,-1){0.2}}

\put(100,127){\line(0,-1){24}}
\put(100,110){\vector(0,-1){0.2}}

\qbezier(100,97.1)(100,90)(105,95) 
\qbezier(103, 100)(110,100)(105,95) 
\put(105,95){\vector(1,1){0.2}}

\put(100.00,130.00){\makebox(0,0)[cc]{$\alpha$}}
\put(100.00,100.00){\makebox(0,0)[cc]{$\epsilon$}}
\put(40.00,100.00){\makebox(0,0)[cc]{$p$}}
\put(40.00,130.00){\makebox(0,0)[cc]{$q$}}
\put(107.00,115.00){\makebox(0,0)[cc]{$0/0$}}
\put(56.00,140.00){\makebox(0,0)[cc]{$1/1$}}
\put(84.00,120.00){\makebox(0,0)[cc]{$1/1$}}
\put(70.00,93){\makebox(0,0)[cc]{$01/10$}}

\put(33.00,115){\makebox(0,0)[cc]{$0/0$}}
\put(110,90){\makebox(0,0)[cc]{$01/01$}}
\put(25.00,145.00){\makebox(0,0)[cc]{$B_4$}}
\end{picture}
\caption{Machine $B_4$.}
\label{at01}
\end{figure}

Let $(q,x,y)\in Q\times A^\omega\times B^\omega$. 
We write $y=q*x$ if $\forall n \in \mathbb{N} \; y[0,n]=q*x[0,n]$ and say machine $V$ {\em transforms} $x$ to $y$.  We refer to words $x$ and $y$ as machines {\em input} and {\em output}, respectively. 

\begin{example} 
\textnormal{
Look at \ref{at01}\ fig. for example of machine $B_4$. }
\end{example}

We might refer to operations $\circ$ and $*$ as machine {\em transition} and {\em output functions}, respectively. 

\begin{definition}
A 3-sorted algebra $V_0=\langle Q,A,B,q_0,\circ ,\ast \rangle$ 
is called  an 
{\em initial Mealy machine} if $\langle Q,A,B,\circ ,\ast \rangle$ is a 
Mealy machine and $q_0\in Q$.
\end{definition}

Suppose that we are given two initial machines \\
$V=\langle Q,A,B;q_0,\circ,*\rangle$ and $V'=\langle Q',A',B';q'_0,\acute\circ,\acute*\rangle$. Schematically it is shown in \ref{z56a}a.\,fig.  

\begin{figure}[h]
\scalebox{.95}{\input{aut28.pic}}
\caption{Serial composition.}
\label{z56a}
\end{figure}

We want to connect the output of machine $V$ to the input of machine $V'$  (shown in \ref{z56a}b.\,fig.). Clearly, in this situation, we have $v=v'$.

Suppose that $B\subseteq A'$, then for the input of the machine $V'$ we always can use the word $v=q_0*u$. Therefore the word $w$ is correctly defined as
\[
w\apz q'_0\acute*(q_0*u).
\]

The symbol $\apz$ is used to make a definition.

\section{Morphism}

We define the morphism $\eta:\{p, q,\alpha\}^+\to \{p, q,\alpha\}^+$
 as follows:
\begin{eqnarray*}
p &\mapsto& pqp\\
q &\mapsto& \alpha\\
\alpha &\mapsto& q
\end{eqnarray*}   

We set 
\begin{eqnarray*}
\eta^0(p) &\apz& p \\
\eta^{\ell+1}(p) &\apz& \eta^\ell(\eta(p))
\end{eqnarray*}

\begin{lemma}\label{l5a.5.1}
$\eta^\ell(p)=\eta^{\ell-1}(p)\delta\eta^{\ell-1}(p)$, 
where
\[
\delta=\begin{cases}
q, & {\rm if} \quad \ell\equiv 1\mod 2,\\
\alpha & {\rm if} \quad \ell\equiv 0\mod 2.
\end{cases}
\]
\end{lemma}

$\Box$ The proof is inductive.
 \begin{eqnarray*}
\eta^1(p) &=& \eta(p)=pqp=\eta^0(p)q\eta^0(p),\\
\eta^2(p) &=& \eta(pqp)=pqp\alpha pqp=\eta^1(p)\alpha \eta^1(p)
\end{eqnarray*}

$\eta^{\ell+1}(p)=\eta(\eta^\ell(p))=\eta(\eta^{\ell-1}(p)\delta'\eta^{\ell-1}(p))=\eta^\ell(p)\eta(\delta')\eta^\ell(p)$.\\
Since
 $\delta'\in\{q,\alpha\}$, it follows that $\eta(\delta')\in\{q,\alpha\}$.

Let $\ell\equiv 0\mod 2$ and $\eta^\ell(p)=\eta^{\ell-1}(p)\alpha\eta^{\ell-1}(p)$, then
\begin{eqnarray*}
\eta^{l+1}(p) &=&\eta(\eta^\ell(p))=\eta(\eta^{\ell-1}(p)\alpha\eta^{\ell-1}(p))=\eta^\ell(p)\eta(\alpha)\eta^\ell(p)=
\eta^\ell(p)q\eta^\ell(p),\\
\eta^{\ell+2}(p) &=&\eta(\eta^\ell(p)q\eta^\ell(p))=\eta^{\ell+1}(p)\alpha\eta^{\ell+1}(p).
\qquad \rule{2mm}{2mm}
\end{eqnarray*}
\medskip

{\bf Further} we are interested exclusively in the machine $B_4$.

{\bf Convention.} We adopt the notational conventions:
\begin{eqnarray*}
\forall v\in\overline{0,1}^+ \quad v^0 &\apz& \lambda \;\wedge\; v^{\ell+1}\apz v^\ell\# v\\
\overline{0,1}^\infty &\apz& \overline{0,1}^*\cup\overline{0,1}^\omega\\
Q &\apz& \{p,q,\alpha, \epsilon\}\\
\forall x\in\overline{0,1}^\infty\; \forall \delta\in Q\quad x\bar \delta &\apz& \delta*x\\
\forall \sigma\in Q^* \quad x\overline{\sigma \delta}  &\apz& (x\bar \sigma)\bar \delta\\
\eta^0(p) &\apz& p\\
\eta^\ell &\apz& \overline{\eta^\ell(p)}
\end{eqnarray*}

\begin{corollary}\label{s5a.5.2}
If $\eta^{\ell+1}(p)=\eta^\ell(p)\delta\eta^\ell(p)$, then
\begin{eqnarray*}
1^\ell01^\omega\bar \delta\eta^\ell &=& 1^\ell001^\omega\eta^\ell,\\
1^{\ell}001^\omega\bar \delta\eta^\ell &=& 1^\ell01^\omega\eta^\ell.
\end{eqnarray*}
and
\[
\delta\circ 1^l00=\epsilon=\delta\circ 1^l01
\]
\end{corollary}

$\Box$ This follows imediately from the fact that
\[
\delta=\left.\begin{cases}
q, & {\rm ja} \quad \ell+1\equiv 1\mod 2,\\
\alpha & {\rm ja} \quad \ell+1\equiv 0\mod 2.
\end{cases}\right\} = 
\begin{cases}
q, & {\rm ja} \quad \ell\equiv 0\mod 2,\\
\alpha & {\rm ja} \quad \ell\equiv 1\mod 2.
\end{cases} \qquad \rule{2mm}{2mm}
\]

\section{Group $\Gamma(B_4)$}

We denote by $\Gamma(B_4)$ the group generated by the set $\{\bar p, \bar q, \bar \alpha, \bar \epsilon\}$, namely,
$\Gamma(B_4)=\langle \bar p, \bar q, \bar \alpha, \bar \epsilon \rangle$. For details see \cite{buls}.
\begin{lemma}\begin{itemize}
\item[\rm (i)]
$1^\omega\eta^\ell=1^\ell01^\omega$, \quad $1^\ell01^\omega\eta^\ell=1^\omega$. 
\item[\rm (ii)] Let $\eta^\ell(p)=p_1p_2\cdots p_m$ and $\eta_j^\ell\apz\overline{p_1p_2\cdots p_j}$, 
\linebreak
$\eta^\ell_0 \apz \mathbb{I} : \overline{0,1}^\infty \to \overline{0,1}^\infty: x\mapsto x$, then 
\begin{eqnarray*}
\overline{0,1}^{\ell+1} &=& \{1^{\ell+1}\eta_j^\ell\,|\,j\in\overline{0,m}\},\\
\overline{0,1}^{\ell+1} &=& \{1^{\ell}0\eta_j^\ell\,|\,j\in\overline{0,m}\}.
\end{eqnarray*}
\item[\rm (iii)] Let $u_j\apz 1^{\ell+1}\eta^\ell_j$,  $v_j\apz 1^\ell0\eta^\ell_j$, then for all indices  $j<m$
\[
p_{j+1}\circ u_j=\epsilon, \quad p_{j+1}\circ v_{j}=\epsilon.
\]
\end{itemize}
\end{lemma}

$\Box$ The proof is inductive. {\bf The induction basis.}
 
 From definitions

\begin{itemize}
\item $\eta^0(p)=p$, \; $\eta^0=\bar p$,
\begin{eqnarray*}
1^\omega\eta^0 &=& 1^\omega \bar p= p*1^\omega=01^\omega,\\
01^\omega\eta^0 &=& 01^\omega\bar p= p*01^\omega=1^\omega.
\end{eqnarray*}

\item $\eta^0_0=\mathbb{I}$, \; $\eta^0_1=\bar p$
\begin{eqnarray*}
\{1\eta^0_0,1\eta^0_1\} &=& \{1,0\}=\overline{0,1},\\
\{0\eta^0_0,0\eta^0_1\} &=& \{0,1\}=\overline{0,1}.
\end{eqnarray*}

\item $p\circ 1=\epsilon, \; p\circ 0=\epsilon$.
\end{itemize}

\begin{itemize}
	\item $\eta^1(p)=\eta(p)=pqp$, \; $\eta^1=\overline{pqp}$.
	\begin{eqnarray*}
1^\omega\eta^1 &=& 1^\omega \overline{pqp}=( p*1^\omega)\overline{qp}=01^\omega\overline{qp}=(q*01^\omega) \bar p\\
&=& 001^\omega\bar p=p*001^\omega=101^\omega,\\
101^\omega\eta^1 &=& 101^\omega\overline{pqp} = (p*101^\omega)\overline{qp}=001^\omega\overline{qp}=(q*001^\omega)\bar p\\
&=& 01^\omega \bar p=p*01^\omega=1^\omega.
\end{eqnarray*}
	
	\item $\eta^1_0=\mathbb{I}, \; \eta^1_1=\bar p, \; \eta^1_2=\overline{pq}, \; \eta^1_3=\overline{pqp}=\eta^1$.
	\begin{eqnarray*}
\{11\eta^1_0,11\eta^1_1,11\eta_2^1,11\eta_3^1\} &=& \{11,01, 00,10\}=\overline{0,1}^2,\\
\{10\eta^1_0,10\eta^1_1,10\eta_2^1,10\eta_3^1\} &=& \{10,00,01,11\}=\overline{0,1}^2.
\end{eqnarray*}

\item 
\begin{align*}
&p\circ 11=\epsilon\circ 1=\epsilon, & q\circ 01=p\circ 1=\epsilon, && p\circ 00=\epsilon\circ0=\epsilon;\\
&p\circ 10=\epsilon\circ 0=\epsilon, & q\circ 00=p\circ 0=\epsilon, && p\circ 01=\epsilon\circ1=\epsilon.	
\end{align*}
\end{itemize}

{\bf The induction step.}
\begin{itemize}
\item
 $\eta^{\ell+1}(p)\underset{\rm{L}\ref{l5a.5.1}}{=}\eta^\ell(p)\delta\eta^\ell(p)=p_1p_2\cdots p_m\delta p_1p_2\cdots p_m$,\\
where $\delta\in\{q,\alpha\}$. 
Hence
\begin{eqnarray*}
1^\omega\eta^{\ell+1}=1^\omega\eta^\ell\bar\delta \eta^\ell=1^\ell01^\omega \bar \delta \eta^\ell \underset{\rm{S}\ref{s5a.5.2}}{=}1^\ell001^\omega\eta^\ell
\end{eqnarray*}
 We have $p_{j+1}\circ v_j=\epsilon$ and $1^\ell01^\omega\eta^\ell=1^\omega$.  Hence 
$1^\ell001^\omega\eta^\ell=1^{\ell+1}01^\omega$.
 Similarly
\[
1^{\ell+1} 01^\omega\eta^{\ell+1}=1^{\ell+1} 01^\omega\eta^\ell\bar\delta\eta^\ell= 1^{\ell} 001^\omega\bar\delta\eta^\ell
\]
because we have $p_{j+1}\circ u_j=\epsilon$ and $1^\omega\eta^\ell=1^\ell01^\omega.$
\[
1^{\ell} 001^\omega\bar\delta\eta^\ell \underset{\rm{S}\ref{s5a.5.2}}{=} 1^{\ell} 01^\omega\eta^\ell=1^\omega
\]

\item
 \[
\eta_j^{\ell+1}=
\begin{cases}
\eta_j^\ell, & {\rm if} \; j\le m,\\
\eta^\ell \bar \delta,  & {\rm if} \; j+m+1,\\
\eta^\ell\bar\delta\eta^\ell_i,  & {\rm if} \; j=m+1+i\wedge i>0.
\end{cases}
\]
Hence
\begin{eqnarray*}
1^{\ell+2}\eta_j^{\ell+1} &=&
\begin{cases}
1^{\ell+2}\eta_j^\ell, & {\rm if} \; j\le m,\\
1^{\ell+2}\eta^\ell \bar \delta,  & {\rm if} \; j+m+1,\\
1^{\ell+2}\eta^\ell\bar\delta\eta^\ell_i,  & {\rm if} \; j=m+1+i\wedge i>0.
\end{cases}
\\
&=&
\begin{cases}
1^{\ell+1}\eta_j^\ell1, & {\rm if} \; j\le m,\\
1^{\ell} 00 ,  & {\rm if} \; j+m+1,\\
1^{\ell}0\eta^\ell_i0,  & {\rm if} \; j=m+1+i\wedge i>0.
\end{cases}
\end{eqnarray*}
We took in consideration that for all indices $j<m$
\[
p_{j+1}\circ u_j=\epsilon, \quad p_{j+1}\circ v_{j}=\epsilon;
\]
furthermore $1^\omega\eta^\ell_m=1^\omega\eta^\ell=101^\omega$ and $1^\ell01^\omega\eta_m^\ell=1^\ell01^\omega\eta^\ell=1^\omega$.
Thus 
\begin{eqnarray*}
1^{\ell+2}\eta_j^{\ell+1} &=&
\begin{cases}
1^{\ell+1}\eta_j^\ell1, & {\rm if} \; j\le m,\\
1^{\ell} 00 ,  & {\rm if} \; j=m+1,\\
1^{\ell}0\eta^\ell_i0,  & {\rm if} \; j=m+1+i\wedge i>0.
\end{cases}\\
&=& 
\begin{cases}
u_j1, & {\rm if} \; j\le m,\\
1^{\ell} 00 ,  & {\rm if} \; j=m+1,\\
v_i0,  & {\rm if} \; j=m+1+i\wedge i>0.
\end{cases} \\
&=&
\begin{cases}
u_j1, & {\rm if} \; j\le m,\\
v_i0,  & {\rm if} \; j=m+1+i\wedge i\ge 0.
\end{cases}
\end{eqnarray*}

Since $\overline{0,1}^{l+2}=\overline{0,1}^{l+1}1\cup\overline{0,1}^{l+1}0$ then we have 
proved that
\[
\overline{0,1}^{\ell+2} = \{1^{\ell+2}\eta_j^{\ell+1}\,|\,j\in\overline{0,2m+1}\}
\]
Similarly 
\begin{eqnarray*}
1^{\ell+1}0\eta_j^{\ell+1} &=&
\begin{cases}
1^{\ell+1}\eta_j^\ell0, & {\rm if} \; j\le m,\\
1^{\ell+2}  ,  & {\rm if} \; j+m+1,\\
1^{\ell+1}\eta^\ell_i1,  & {\rm if} \; j=m+1+i\wedge i>0.
\end{cases}\\
&=& 
\begin{cases}
v_j0, & {\rm if} \; j\le m,\\
1^{\ell+2}  ,  & {\rm if} \; j+m+1,\\
u_i1,  & {\rm if} \; j=m+1+i\wedge i>0.
\end{cases}\\
&=&
\begin{cases}
v_j0, & {\rm if} \; j\le m,\\
u_i1,  & {\rm if} \; j=m+1+i\wedge i\ge 0.
\end{cases}
\end{eqnarray*}
Therefore
\[
\overline{0,1}^{\ell+2} = \{1^{\ell+1}0\eta_j^{\ell+1}\,|\,j\in\overline{0,2m+1}\}
\]

\item Let $\dot u_j\apz 1^{\ell+2}\eta_j^{\ell+1}, \; \dot v_j=1^{l+1}0\eta_j^{\ell+1}$. We must prove that for all $j<2m+1$
\[
q_{j+1}\circ \dot u_j=\epsilon, \quad q_{j+1}\circ \dot v_j=\epsilon, 
\]
where $\eta^{\ell+1}(p)=q_1q_2\cdots q_{2m+1}$.

We know
\begin{eqnarray*}
q_1q_2\cdots q_m &=& p_1p_2\cdots p_m,\\
q_{m+1} &=& \delta, \\
q_{m+2}q_{m+3}\cdots q_{2m+1} &=& p_1p_2\cdots p_m.
\end{eqnarray*}
\[
\dot u_j=\begin{cases}
u_j1, & {\rm if} \; j\le m,\\  
v_i0, & {\rm if} \; j=m+i.
\end{cases} \quad
\dot v_j=\begin{cases}
v_j0, & {\rm if} \; j\le m,\\
u_i1, & {\rm if} \; j=m+i.
\end{cases}
\]
In particular
\begin{align*}
	&\dot u_m=1^{\ell+2}\eta_m^{\ell+1}=1^{\ell+2}\overline{p_1p_2\cdots p_m}=1^{\ell+2}\eta^\ell=1^\ell01\\
	&\dot v_m=1^{\ell+1}0\eta_m^{\ell+1}=1^{\ell+1}0\overline{p_1p_2\cdots p_m}=1^{\ell+1}0\eta^\ell=1^{\ell}00
\end{align*}
Thus
\begin{align*}
 \dot u_{j+1}&=p_{j+1}* \dot u_j, & &{\rm if} \quad j\in\overline{0,m-1},\\
 \dot u_{m+1}&=\delta* \dot u_m,  \\
 \dot u_{j+1+m}&=p_j* \dot u_{j+m}, & &{\rm if} \quad j\in\overline{1,m}.
\end{align*}
Subsequently
\begin{eqnarray*}
 &&q_{j+1}\circ \dot u_j=\\
&&=
\begin{cases}
p_{j+1}\circ u_j1=p_{j+1}\circ u_j\circ 1=\epsilon\circ 1=\epsilon, &{\rm if} \; j\in\overline{0,m-1},\\
\delta\circ \dot u_m=\delta\circ 1^l01 \underset{\rm{S}\ref{s5a.5.2}}{=} \epsilon, &{\rm if} \; j=m,\\
p_{i+1}\circ v_i0=p_{i+1}\circ v_i\circ 0=\epsilon\circ 0=\epsilon, &{\rm if} \; i\in\overline{0,m-1}\wedge j=m+1+i.
\end{cases}
\end{eqnarray*}
Similarly 
\begin{eqnarray*}
&& q_{j+1}\circ \dot v_j=\\
&&=\begin{cases}
p_{j+1}\circ v_j0   =p_{j+1}\circ v_j\circ 0=\epsilon\circ 0=\epsilon, &{\rm if} \; j\in\overline{0,m-1},\\
\delta\circ \dot v_m=\delta\circ 1^l00 \underset{\rm{S}\ref{s5a.5.2}}{=} \epsilon, &{\rm if} \; j=m,\\
p_{i+1}\circ u_i1=p_{i+1}\circ u_i\circ 1=\epsilon\circ 1=\epsilon, &{\rm if} \; i\in\overline{0,m-1}\wedge j=m+1+i.
\end{cases}
\end{eqnarray*}
\end{itemize}
This completes the induction. \rule{2mm}{2mm}

\begin{corollary} \label{s5a.5.4}
Group $\Gamma(B_4)$ is infinite.
\end{corollary}

$\Box$ Since $1^\omega\eta^\ell=1^\ell01^\omega$ then all elements $\eta^\ell$ of $\Gamma(B_4)$ are distinct.~\rule{2mm}{2mm}

\section{$\Gamma(B_4)$ is not periodic.}

\begin{definition}
A group is called periodic if every element of the group has finite order.
\end{definition}

\begin{lemma}
$\langle \bar\alpha, \bar q\rangle \cong \mathbb{Z}_2\times \mathbb{Z}_2$
\end{lemma}

$\Box$ (i) $\overline{\alpha^2}=\mathbb{I}=\overline{q^2}$.
Let $x=1^\ell0x_1x_2x_3\cdots$ then
\[
x\bar \alpha \bar\alpha = \left.\begin{cases}
1^\ell 0x_1x_2x_3\cdots, & {\rm if} \; \ell \equiv 0\mod 2\\
1^\ell0\tilde x_1x_2x_3\cdots, & {\rm if} \; \ell \equiv 1\mod 2
\end{cases}\right\} \bar\alpha= 1^\ell0x_1x_2x_3\cdots=x
\]
\[
x\bar q \bar q = \left.\begin{cases}
1^\ell 0x_1x_2x_3\cdots, & {\rm if} \; \ell \equiv 1\mod 2\\
1^\ell0\tilde x_1x_2x_3\cdots, & {\rm if} \; \ell \equiv 0\mod 2
\end{cases}\right\} \bar q= 1^\ell0x_1x_2x_3\cdots=x
\]
Here
\[
\tilde x_1\apz\begin{cases}
0, & {\rm if} \; x_1=1;\\
1, & {\rm if} \; x_1=0.
\end{cases}
\]
\[
x\bar \alpha \bar q = \left.\begin{cases}
1^\ell 0x_1x_2x_3\cdots, & {\rm if} \; \ell \equiv 0\mod 2\\
1^\ell0\tilde  x_1x_2x_3\cdots, & {\rm if} \; \ell \equiv 1\mod 2
\end{cases}\right\} \bar q = 1^\ell0 \tilde  x_1x_2x_3\cdots
\]
\[
x\bar q \bar \alpha = \left.\begin{cases}
1^\ell 0x_1x_2x_3\cdots, & {\rm if} \; \ell \equiv 1\mod 2\\
1^\ell0\tilde  x_1x_2x_3\cdots, & {\rm if} \; \ell \equiv 0\mod 2
\end{cases}\right\} \bar \alpha= 1^\ell0\tilde  x_1x_2x_3\cdots
\]
Thus $\overline{\alpha q}=\overline{q \alpha}$.
Hence $\langle \bar \alpha, \bar q \rangle=\{\mathbb{I}, \bar \alpha, \bar q, \overline{\alpha q}\}$
because words from $\{\alpha, q\}^3$ do not generate new elements. For example
$\overline{\alpha q\alpha}=\bar \alpha\overline{q\alpha}=\bar \alpha\overline{\alpha q}=\overline{\alpha \alpha q}=\mathbb{I}\bar q=\bar q$.

There are only 2 groups (up to isomorphism) of order 4.  The group $\langle \bar \alpha, \bar q \rangle$ is not the cyclic group. Therefore $\langle \bar\alpha, \bar q\rangle \cong \mathbb{Z}_2\times \mathbb{Z}_2$.     
\rule{2mm}{2mm}
\medskip

We can assume that every element $g$ of $\Gamma(B_4)$ is represented as word 
\[
w=sa_1pa_2p\cdots pa_n\sigma
\]
where $a_i\in \{q,\alpha, \beta\}$ and $s,\sigma\in\{\lambda,p\}$.
Here $g=\bar w$ and $\bar\beta\apz\overline{\alpha q}$. We took in consideration that the order of elements 
$\bar p, \bar q, \bar \alpha, \bar\beta$ is 2; 
furthermore the elements  $\bar q,\bar\alpha,\bar\beta$ commute with each other and 
\[
\overline{\alpha q}=\bar\beta, \; \overline{q\beta}=\overline{q\alpha q}=\overline{\alpha qq}=\bar\alpha, \; \overline{\beta\alpha}=\overline{\alpha q\alpha}=\overline{\alpha\alpha q}=\bar q.
\]

\begin{definition}
Let $G$ be a group and let $a,b\in G$. Then $a$ is conjugate to $b$ if there is a $h\in G$ such that $b= hah^{-1}$.
\end{definition}

Let 
\[
S(a)\apz\{b\,|\,\exists h\in G \quad b=hah^{-1}\}
\]
denotes the conjugacy class of the element $a$. 

\begin{lemma}\label{l5a.5.7}
Let $G$ be a group and let $a,b\in G$. If $a$ has the finite order $o(a)=n$ and $b\in S(a)$ then $o(a)\le n$.
\end{lemma}

$\Box$ $b^n=(gag^{-1})^n=gag^{-1}gag^{-1}\cdots gag^{-1}=ga^ng^{-1}=geg^{-1}=e$.
Here $e$ is the neutral element of the group $G$.
\rule{2mm}{2mm}

\begin{lemma}
\begin{eqnarray*}
o(\bar p) &=& o(\bar q)=o(\bar \alpha)=o(\overline{\alpha q})=2,\\
o(\overline{qp})  &\le& o(\overline{pq}) = 8,\\
o(\overline{\alpha p}) &\le& o(\overline{p\alpha})=4
\end{eqnarray*}
\end{lemma}

$\Box$
(i)  Let $x=x_0x_1\cdots x_n\cdots\in\overline{0,1}^\omega$ and $y=x_1x_2\cdots x_n\cdots$ then
\begin{eqnarray*}
001^\ell 0x \overline{pq} &=& 101^\ell0x\bar q=101^\ell0x,\\
101^\ell0x\overline{pq} &=& 001^\ell 0x\bar q=01^{\ell+1}0x
\end{eqnarray*}

 Case 1:  $\ell \equiv 0\mod 2$
\begin{eqnarray*}
01^{\ell+1}0x\overline{pq} &=& 1^{\ell+2}0x\bar q= 1^{\ell+2}0\tilde  x_0y,\\
1^{\ell+2}0\tilde  x_0y\overline{pq} &=& 01^{\ell+1}0\tilde  x_0y \bar q=001^\ell0\tilde  x_0y,\\
001^\ell0\tilde  x_0y \overline{pq} &=& 101^\ell0\tilde  x_0y\bar q=101^{\ell}0 \tilde  x_0y,\\
101^\ell0 \tilde  x_0y \overline{pq} &=& 001^\ell0 \tilde  x_0 y \bar q=01^{\ell+1}0 \tilde  x_0 y,\\
01^{\ell+1}0 \tilde  x_0 \overline{pq} &=& 1^{\ell+2}0 \tilde  x_0 y \bar q= 1^{\ell+2}0x_0y=1^{\ell+2}0x,\\
1^{\ell+2}0x \overline{pq} &=& 01^{\ell+1}0x\bar q=001^\ell0x
\end{eqnarray*}
Hence if \[
z\in\{001^\ell0x, 101^\ell0x, 011^\ell0x,111^\ell0x\}
\]
then $z(\overline{pq})^8=z$.
\medskip

Case 2: $\ell\equiv 1\mod 2$
\begin{eqnarray*}
001^\ell 0x \overline{pq} &=& 101^\ell0x\bar q=101^\ell0x,\\
101^\ell0x\overline{pq} &=& 001^\ell 0x\bar q=01^{\ell+1}0x,\\
01^{\ell+1}0x\overline{pq} &=& 1^{\ell+2}0x\bar q= 1^{\ell+2}0 x,\\
1^{\ell+2}0 x\overline{pq} &=& 01^{\ell+1}0 x \bar q=001^\ell0 x
\end{eqnarray*}
Hence if \[
z\in\{001^\ell0x, 101^\ell0x, 011^\ell0x,111^\ell0x\}
\]
then $z(\overline{pq})^4=z$.

What happens with word $001^\omega$?
\begin{eqnarray*}
001^\omega\overline{pq} &=& 101^\omega\bar q=101^\omega,\\
101^\omega\overline{pq} &=& 001^\omega\bar q=01^\omega,\\
01^\omega\overline{pq} &=& 1^\omega\bar q= 1^\omega,\\
1^\omega\overline{pq} &=& 01^\omega \bar q=001^\omega
\end{eqnarray*}
Hence if 
\[
z\in\{001^\omega, 101^\omega, 01^\omega, 1^\omega\}
\]
then $z(\overline{pq})^4=z$.

This completes the proof for $(\overline{pq})^8=\mathbb{I}$.

Now
\[
(\overline{qp})^8=(\bar p)^2(\overline{qp})^8=\bar p((\overline{pq})^8)\bar p.
\]
Thus (see Lemma \ref{l5a.5.7}) $o(\overline{qp})\le 8$.
\medskip

(ii) It seems that $(\overline{p\alpha})^8=\mathbb{I}$ but we need the proof.

Case 1:
\begin{eqnarray*}
001x\overline{p\alpha} &=& 101x \bar\alpha=100x,\\
100x\overline{p\alpha} &=& 000x\bar \alpha=000x,\\
000x\overline{p\alpha} &=& 100x\bar \alpha= 101x,\\
101x\overline{p\alpha} &=& 001x\bar\alpha=001x
\end{eqnarray*}
Hence if 
\[
z\in\{00x, 10x\}
\]
then $z(\overline{p\alpha})^4=z$.
\medskip

Case 2:
\begin{eqnarray*}
1^{\ell+2} 0x \overline{p\alpha} &=& 01^{\ell+1}0x\bar \alpha=01^{\ell+1}0x\\
\end{eqnarray*}
\medskip

Case 2a: $\ell \equiv 1\mod 2$
\begin{eqnarray*}
01^{\ell+1}0x \overline{p\alpha} &=& 1^{\ell+2}0x\bar \alpha=1^{\ell+2}0\tilde  x_0y\\
1^{\ell+2}0\tilde  x_0y &=& 01^{\ell+1}0\tilde x_0y\bar \alpha=01^{\ell+1}0\tilde x_0y,\\
01^{\ell+1}0\tilde x_0y \overline{p\alpha} &=& 1^{\ell+2}0\tilde x_0y=1^{\ell+2}0x_0y=1^{\ell+2}0x
\end{eqnarray*}
Hence if 
\[
z\in\{011^\ell0x, 111^\ell0x\}
\]
then $z(\overline{p\alpha})^4=z$.
\medskip

Case 2b:  $\ell\equiv 0\mod 2$
\begin{eqnarray*}
01^{\ell+1}0x \overline{p\alpha} &=& 1^{\ell+2}0x\bar \alpha=1^{\ell+2}0x
\end{eqnarray*}
Hence if 
\[
z\in\{011^\ell0x, 111^\ell0x\}
\]
then $z(\overline{p\alpha})^2=z$.
\medskip

What happens with word $01^\omega$? 
\begin{eqnarray*}
01^\omega \overline{p\alpha} &=& 1^\omega\bar \alpha=1^\omega,\\
1^\omega \overline{p\alpha} &=& 01^\omega\bar \alpha =01^\omega
\end{eqnarray*}
Hence if 
\[
z\in\{01^\omega, 1^\omega\}
\]
then $z(\overline{p\alpha})^2=z$.

This completes the proof for $(\overline{p\alpha})^4=\mathbb{I}$.

Now
\[
(\overline{\alpha p})^4=(\bar p)^2(\overline{\alpha p})^4=\bar p((\overline{p\alpha})^4)\bar p.
\]
Thus (see Lemma \ref{l5a.5.7}) $o(\overline{\alpha q})\le 4$.
\rule{2mm}{2mm}

\begin{lemma}\label{l5a.5.9}
Let $\xi\apz \overline{p\alpha q}$. If

{\rm (i)} \quad $1^\omega\xi^k=u_kx_k$,

{\rm (ii)} \quad $|u_k|=n$,

{\rm (iii)} \quad $\ell < 2^{n}$,

then

{\rm (i)} \quad $0 \smallsetminus u_\ell$,

{\rm (ii)} \quad $\overline{0,1}^n=\{u_k\,|\,k\in\overline{1,2^n}\}$,

{\rm (iii)} \quad $x_\ell=1^\omega$, ja $\ell<2^{n-1}$,

{\rm (iv)} \quad $x_\ell=01^\omega$, ja $2^{n-1}\le \ell<2^n$,

{\rm (v)} \quad $1^\omega \xi^{2^n}=1^n0^21^\omega$,

\end{lemma}

$\Box$ The proof is inductive. {\bf The induction basis.}
\begin{eqnarray*}
1^\omega \xi &=& 00111^\omega\\
1^\omega \xi^2 &=& 10011^\omega\\
1^\omega \xi^3 &=& 01011^\omega\\
1^\omega \xi^4 &=& 11001^\omega\\
1^\omega \xi^5 &=& 00001^\omega\\
1^\omega \xi^6 &=& 10101^\omega\\
1^\omega \xi^7 &=& 01101^\omega\\
1^\omega \xi^8 &=& 111001^\omega
\end{eqnarray*}

{\bf The induction step} for
 $n\ge 3$.

 Let  $1^\omega \xi^k=v_ky_k$ where  $|v_k|=n+1$. We know $1^\omega \xi^k=u_kx_k$. Therefore
$v_k=u_ka_k$ where $a_k\in\overline{0,1}$.

(i) If $\ell<2^{n}$ then $0 \smallsetminus u_\ell$. Therefore $0 \smallsetminus u_\ell a_\ell=v_\ell$.

If $\ell=2^{n}$ then $1^\omega \xi^{2^n}=1^n0^21^\omega$, $v_\ell=1^n0$. Hence $0 \smallsetminus u_\ell$.

If $2^n<\ell<2^{n+1}$ then $\ell=2^n+t$ where $0<t<2^n$. Look
\[
v_\ell y_\ell=1^\omega\xi^\ell=(1^\omega\xi^{2^n})\xi^t=1^n0^21^\omega\xi^t=(1^n\xi^t)a_\ell y_\ell=u_ta_\ell y_\ell
\]
Since $0 \smallsetminus u_t$ and $u_ta_\ell=v_\ell$ then $0 \smallsetminus v_\ell$.

(ii) Let $\kappa=2^{n-1}$ then $1^\omega\xi^\kappa=1^\omega\xi^{2^{n-1}}=1^{n-1}001^\omega$.
Hence $u_\kappa=1^{n-1}0$. It is possible only then if $u_{\kappa-1}=01^{n-2}0$ because $u_{\kappa-1}\xi=u_\kappa$.
Now look! How does  machine $B_4$ work? We can deduce: if $u_\ell\ne u_{\kappa-1}$ then $u_\ell a\xi=u_{\ell+1}a$ for all
$a\in\overline{0,1}$.
\begin{itemize}
\item If the first letter of $u_\ell$ is 1 then the map $\bar p$ transforms 1 to 0 but the map $\overline{\alpha q}$ now  transforms only the secon letter. Thus $u_\ell a\xi=u_{\ell+1}a$.
\item If the first letter of $u_\ell$ is 0 then the map $\bar p$ transforms 0 to 1. So we have a new word $v$. 
\begin{itemize}
\item[$\rhd$] If $1v=1^n$ then $1^na\overline{\alpha q}=1^na$. Thus $u_\ell a\xi=u_{\ell+1}a$.
\item[$\rhd$]  Let $v=v_1v_2\cdots v_n$ and  $0\smallsetminus v$ then $1v\ne 1^{n-1}0$ because $u_\ell \ne u_{\kappa-1}$. It means the first occurrance  $v_i$ of 0 in $v$ is not $v_n$. We can deduce: the map $\overline{\alpha q}$ transforms only the letter $v_{i+1}$. Thus $u_\ell a\xi=u_{\ell+1}a$. 
\end{itemize}
\end{itemize}
We have $\overline{0,1}^n=\{u_k\,|\,k\in\overline{1,2^n}\}$ therefore in the sequence
\[
u_1,u_2,\dots, u_{2^n}
\]
there is only one word equals  $01^{n-2}0=u_{\kappa-1}$. Hence
\[
\forall \ell<\kappa-1\; u_\ell\ne u_{\kappa-1}.
\]
Thus  $\forall\ell\le\kappa-1\; v_\ell=u_\ell1$ because $v_1=001^{n-1}=u_11$ and $v_\ell=v_{\ell-1}\xi=u_{\ell-1}1\xi=u_\ell1$.

We take in consideration that 
\[
v_\kappa=u_{\kappa-1}1\xi=01^{n-2}01\xi=1^{n-1}00=u_\kappa0.
\]
Besides there is no any element $u_\ell$ equals  $u_{\kappa-1}$ in the sequence $u_\kappa, u_{\kappa+1},\ldots, u_{2^n}$. Consequently 
\[
\forall \ell \ge \kappa \; (\ell \le 2^n \Rightarrow v_\ell=u_\ell 0).
\]

We know $v_\ell=u_ta_\ell$ for $\ell=2^n+t$ where $0<t<2^n$.  Thus if $t<\kappa-1$ then 
$v_{\ell+1}=v_\ell\xi= u_t0\xi=u_{t+1}0$. If $t=\kappa-1$ then 
\[
v_{2^n+\kappa}=v_{2^n+\kappa-1}\xi=u_{\kappa-1}0\xi=01^{n-2}00\xi=1^{n-1}01=u_\kappa1.
\]
Hence for all $\ell=2^n+t$ we have $v_\ell=u_t1$ where $\kappa\le t\le 2^n$.

We have got a list
\begin{align*}
	& v_1=u_11, v_2=u_21, \ldots, v_{\kappa-1}=u_{\kappa-1}1,\\
	& v_\kappa=u_\kappa0, v_{\kappa+1}=u_{\kappa+1}0, \ldots, v_{2^n}=u_{2^n}0,\\
	& v_{2^n+1}=u_10, v_{2^n+2}=u_20, \ldots, v_{2^n+\kappa-1}=u_{\kappa-1}0,\\
	& v_{2^n+\kappa}=u_\kappa1, v_{2^n+\kappa+1}=u_{\kappa+1}1, \ldots, v_{2^{n+1}}=v_{2^n+2^n}=u_{2^n}1.
\end{align*}
Thus $\overline{0,1}^{n+1}=\{v_k\,|\,k\in\overline{1,2^{n+1}}\}$.

(iii) We take in consideration that  $x_k=a_ky_k$.

 If $\ell<2^{n-1}$ then  $x_\ell=1^\omega$ therefore $y_\ell=1^\omega$.

If $2^{n-1}\le \ell<2^n$ then $x_\ell=01^\omega$ therefore $y_\ell=1^\omega$.

(iv)  Let $m\apz 2^n$ then  $1^\omega\xi^m=1^\omega \xi^{2^n}=1^n0^21^\omega$ thence
$v_m=1^{n}0$ and $y_m=01^\omega$. It is possible only if $v_{m-1}=01^{n-1}0$ because $v_{m-1}\xi=v_m$.

Now look! How does  machine $B_4$ work? We can deduce: if 
$v_\ell\ne v_{m-1}$ then $v_\ell a\xi=v_{\ell+1}a$ for all
$a\in\overline{0,1}$.

We are interested in fact: when $v_\ell x\xi=v_{\ell+1}x$ for all
$x\in\overline{0,1}^\omega$?
\begin{itemize}
\item  If the first letter of $v_\ell$ is 1 then the map $\bar p$ transforms 1 to 0 but the map $\overline{\alpha q}$ now  transforms only the secon letter. Thus $v_\ell x\xi=v_{\ell+1}x$.
\item If the first letter of $v_\ell$ is 0 then the map $\bar p$ transforms 0 to 1. So we have a new word $v$. 
\begin{itemize}
\item[$\rhd$] If $1v=1^{n+1}$ then $1^{n+1}a\overline{\alpha q}=1^{n+1}a$, but if $0\smallsetminus x$ then $1^{n+1}x\overline{\alpha q}=1^{n+1}x'$ where $x'\ne x$ nevertheless. 

\item[$\rhd$]  Let $v=v_1v_2\cdots v_{n+1}$ and  $0\smallsetminus v$ then $1v\ne 1^{n}0$ because $v_\ell \ne v_{m-1}$. 
It means the first occurrance  $v_i$ of 0 in $v$ is not $v_{n+1}$. We can deduce: the map $\overline{\alpha q}$ transforms only the letter $v_{i+1}$. Thus $v_\ell x\xi=v_{\ell+1}x$. 
\end{itemize}
\end{itemize}
We have $\overline{0,1}^{n+1}=\{v_k\,|\,k\in\overline{1,2^{n+1}}\}$ therefore in the sequence
\[
v_1,v_2,\dots, v_{2^{n+1}}
\]
there is only one word equals  $01^{n-1}0=v_{m-1}$. 
This means there is no any element $v_\ell$ equals  $v_{m-1}$ in the sequence
\[
v_m,v_{m+1},\dots, v_{2^{n+1}}.
\]
Thus if $\ell\ge 2^n$ and $v_\ell\ne 01^n$ then $v_\ell y_\ell=v_{\ell+1}y_\ell$.

We have $u_{2^n-1}=01^{n-1}$ because $u_{2^n}=1^n$. Hence $v_{2^{n+1}-1}=u_{2^n-1}1=01^n$. Thus if $2^n\le \ell < 2^{n+1} $ then  $y_\ell=01^\omega$.

(v) We konow $1^\omega\xi^{2^{n+1}-1}=v_{2^{n+1}-1}y_{2^{n+1}-1}=01^n01^\omega$ therefore
\linebreak
$1^\omega\xi^{2^{n+1}}=01^n01^\omega\xi=1^{n+1}0^21^\omega$.

This completes the induction. 
\rule{2mm}{2mm}

\begin{corollary} \label{s5a.5.4}
Group $\Gamma(B_4)$ is not periodic.
\end{corollary}

$\Box$ Since $1^\omega \xi^{2^n}=1^n0^21^\omega$ then all elements $\xi^\ell$ of $\Gamma(B_4)$ are distinct.
Thus $o(\xi)=\aleph_0$.~\rule{2mm}{2mm}

\section{Dense orbit}
\begin{definition}  Let $u,v\in A^{\infty}\apz A^*\cup A^\omega$. The mapping $d: A^{\infty} \times A^{\infty} \to {\mathbb R}$ is called a {\it metric} or {\it prefix metric} in the set
 $A^{\infty}$ if $$d(u,v)= \left \{\begin{array}{lr}2^{-m}, & u\ne v,\\0, & u=v,\end{array}\right .$$where$$m=\max\{\,|w|\,|\,w\in\mathrm{Pref}(u)\cap\mathrm{Pref}(v)\}.$$
\end{definition}

\begin{definition}
 Let $X,Y\subseteq W$ and $X\subseteq Y$. Then $X$ is {\it dense }  in $Y$ if for each point $y\in Y$
  and each $\varepsilon >0$, there exists $x\in X$ such that $d(x,y)<\varepsilon$.
\end{definition}

The set $\mathcal{O}(x)\apz\{y\,|\,\exists k\; y=x\xi^k\}$ is called the {\em orbit} of $x\in \overline{0,1}^\omega$.
\begin{proposition}
The orbit of every element $x\in \overline{0,1}^\omega$ is dense in $\overline{0,1}^\omega$.
\end{proposition}

$\Box$ At first we are interested in particular case, namely, $x=1^\omega$. Let $\varepsilon>0$ then we can choose $m$ so large that $2^{-m}<\varepsilon$. Let $x\in\overline{0,1}^\omega$, $u\in{\rm Pref}(x)$ and $|u|=m$, in other words, $u\in\overline{0,1}^m$. Now we take into consideration Lemma \ref{l5a.5.9}:
\[
\exists k\; (1^\omega\xi^k=u_k x_k \wedge u_k=u).
\]
Hence $d(x,u_k x_k)\le 2^{-m}<\varepsilon$.

Let $y\in\overline{0,1}^\omega$, $v\in{\rm Pref}(y)$ and $|v|=m$. 
Let $x\in\overline{0,1}^\omega$, $u\in{\rm Pref}(x)$ and $|u|=|v|$. Now we can deduce from Lemma \ref{l5a.5.9}:
 $\exists k\;  u\xi^k=v$. Thus $v\in{\rm Pref}(x\xi^k)$. Hence $d(x\xi^k,y)\le 2^{-m}<\varepsilon$.
So   $\mathcal{O}(x)$ is dense in $\overline{0,1}^\omega$.
\rule{2mm}{2mm}

\section{Topological transitivity}

\begin{definition}   The function $f: X\to X$ is called
 {topologically transitive} on $X$ if
 \begin{eqnarray*}\forall x,y \in X  \;\forall \varepsilon >0 \; \exists z\in X\; \exists n\in {\mathbb N}  \\
  d(x,z)<\varepsilon \; \wedge \; d (y,f^n (z))<\varepsilon.
\end{eqnarray*} 
\end{definition}

\begin{corollary}
The map $\xi: \overline{0,1}^\omega \to \overline{0,1}^\omega$ is topologically transitive on $\overline{0,1}^\omega$.
\end{corollary}

$\Box$
Let $x,y\in\overline{0,1}^\omega$. We can choose $x$ as a word $z$. Since orbit $\mathcal{O}(x)$ is dense in $\overline{0,1}^\omega$ then for every $\varepsilon>0$ exists  $n$ such that
$d(x\xi^n,y)<\varepsilon$.
\rule{2mm}{2mm}

\section{Sensitivity}

\begin{definition}
 The function $f: X\to X$ { exhibits sensitive dependence on initial conditions} if 
\begin{eqnarray*}\exists \delta >0 \;\forall x \in X \;\forall \varepsilon >0 \;\exists y \in X \;\exists 
n\in{\mathbb N}  \\
 d(x,y)<\varepsilon \; \wedge\;  d(f^n(x), f^n(y))>\delta.\end{eqnarray*}
\end{definition}

\begin{definition} A total mapping $f:A^*\to B^*$ is called a sequential function if 
\begin{itemize}
\item[\textnormal{(i)}] $\forall u\in A^*\;|u|=|f(u)|;$
\item[\textnormal{(ii)}] $u\in \textnormal{Pref}(v)\Rightarrow f(u)\in\textnormal{Pref}(f(v))$.
\end{itemize}
\end{definition}

\begin{corollary} For all sequential functions, we have that if
\[u\in\textnormal{Pref}(v)\cap\textnormal{Pref}(w),\]
then
\[f(u)\in\textnormal{Pref}(f(v))\cap\textnormal{Pref}(f(w)).\] 
\end{corollary}

It states that if words $u$ and $v$ have matching prefixes of length $k$, then words $f(u)$ and $f(v)$ have matching prefixes of length $k$.

$\Box$ Suppose that $u\in\textnormal{Pref}(v)\cap\textnormal{Pref}(w)$, then accordingly with the definition of sequential function $f(u)\in\textnormal{Pref}(f(v))$ and $f(u)\in\textnormal{Pref}(f(w))$.~\rule{2mm}{2mm}

\begin{proposition}
If $f: A^\omega \to A^\omega$ is a sequential function then $f$ does not  exhibit sensitive dependence on initial conditions.
\end{proposition}

$\Box$ Let $d(x, y)<\varepsilon$ then exists $m$ such that $d(x, y)= 2^{-m}\le\varepsilon$.
This means that 
$x=ux'$ un $y=uy'$ for some $x',y'\in A^\omega$ where $|u|=m$. Since $f(x)= f(ux') = f(u) x''$ and $f(y)=f(uy')=f(u) y''$
for some $x'',y''\in A^\omega$ then
\[
\forall n \; d(x,y)\ge d(f^n(x),f^n(y)).
\]
Thus $\forall n \; d(f^n(x),f^n(y))<\delta$  for all $\varepsilon<\delta$.
\rule{2mm}{2mm}

\begin{corollary}
The map $\xi: \overline{0,1}^\omega \to \overline{0,1}^\omega$  does not  exhibit sensitive dependence on initial conditions.
\end{corollary}

\end{document}